\documentclass[12pt]{article}
\usepackage{graphicx}
\usepackage{amsmath}
\usepackage[dvips]{epsfig}
\usepackage{amssymb}

\makeatletter
\renewcommand{\section}{\@startsection
  {section}%
  {2}%
  {0mm}%
  {\baselineskip}%
  {0.5 \baselineskip}%
  {\centering}}
\makeatother
\begin{document}

\title { On the weighted $q$-Bernoulli numbers and  polynomials  }
\author{  T.  Kim$$ \\
$^1$ Division of General Education-Mathematics,\\
 Kwangwoon University, Seoul 139-701,  Korea \\
                      }

\date{}
\maketitle

 {\footnotesize {\bf Abstract}\hspace{1mm}
In this paper we discuss  new concept of the  $q$-extension of
Bernoulli numbers and polynomials with weight $\alpha$. From these
$q$-Bernoulli numbers with weight $\alpha$, we establish some
interesting identities and relations.



\begin{center} {\bf 1. Introduction} \end{center}

Let $p$ be a fixed  odd prime.  Throughout this paper
$\mathbb{Z}_p$, $\Bbb Q_p$, $\mathbb{C}$ and  $\mathbb{C}_p$
will, respectively, denote the ring of $p$-adic rational integers,
the field of  $p$-adic rational numbers,  the complex number field
and the completion of algebraic closure of $\Bbb Q_p$.
 Let $\nu_p$ be the normalized exponential valuation of
$\mathbb{C}_p$ with $|p|_p=p^{-\nu_p(p)}=p^{-1}$. When one talks
of $q$-extension, $q$ is variously considered as an indeterminate,
a complex number $q\in \mathbb{C},$ or $p$-adic number $q\in\Bbb
C_p .$ If $q\in \Bbb C,$  one normally assume $|q|<1.$ If $q\in
\mathbb{C}_p,$ then we assume  $|q-1|_p<p^{-\frac{1}{p-1}},$ so
that $q^x=\exp(x \log q)$ for $|x|_p\leq 1.$  In this paper we use
the following notation:
$$[x]_q =\frac{1-q^x}{1-q}, \mbox{ (see [1-16])}.$$
Let $UD(\mathbb{Z}_p)$ be the space of uniformly differentiable
functions on $\mathbb{Z}_p.$
 For
$ f \in UD(\mathbb{Z}_p)$,
 let us start with the expressions
$$ \dfrac{1}{[p^N]_q} \sum_{ 0 \leq x < p^N} q^x f(x)=\sum_{ 0 \leq x < p^N} f(x) \mu_q(x + p^N \Bbb Z_p),$$
representing $p$-adic $q$-analogue of Riemann sums for $f$ (see
[3-7]).

The integral of $f$ on $\Bbb Z_p$ will be defined as the limit( $N
\rightarrow \infty$) of these sums, when it exists.
 The  $p$-adic $q$-integral of function $f$($ \in UD(\mathbb{Z}_p)$) is  defined by
$$I_q(f)=\int_{\mathbb{Z}_p}f(x) d\mu_{q}(x)=
\lim_{N \to \infty}\dfrac{1}{[p^N]_q}\sum_{x=0}^{ p^N-1}f(x)q^x,
\text{ (see [7])}. \eqno(1)
$$
From (1), we note that
$$\big| \int_{\mathbb{Z}_p}f(x) d\mu_{q}(x) \big| \leq p \| f \|_1,  \mbox{ (see
[3-7])}, $$ where $$  \| f \|_1= \sup \left \{ | f(0)|_p, \quad
\sup_{x \neq y} \big| \dfrac{f(x)-f(y)}{x-y} \big|_p \right \}.$$

In [2], Carlitz defined  a set of numbers $\zeta_k=\zeta_k(q)$
inductively by
$$\zeta_{0, q}=1,  \quad  ( q \zeta+1)^k -  \zeta_{k, q} =
\left \{\begin{array}{ll}
1, & \mbox{ if } k=1, \\
0, &  \mbox{ if } k>1,
\end{array} \right. \eqno(2)
$$ with the usual convention of replacing $\zeta^k$ by $\zeta_{k, q}$.
These numbers are the $q$-extension of ordinary Bernoulli numbers.
But they do not remain finite when $q=1$. So,  Carlitz modified
Eq.(2) as follows:
$$\beta_{0, q}=1,  \quad  q ( q \beta+1)^k -  \beta_{k, q} =
\left \{\begin{array}{ll}
1, & \mbox{ if } k=1, \\
0, &  \mbox{ if } k>1,
\end{array} \right.
$$ with the usual convention of replacing $\beta^k$ by $\beta_{k, q}$(see [2]).
These numbers $\beta_{k, q}$  are called the $n$-th Carlitz's
$q$-Bernoulli numbers.

In [1], Carlitz also considered the extension of Carlitz's
$q$-Bernoulli numbers as follows:
$$\beta_{0, q}^{h}=\dfrac{h}{[h]_q},  \quad  q^{h} ( q \beta^{h}+1)^k -  \beta_{k, q}^{h} =
\left \{\begin{array}{ll}
1, & \mbox{ if } k=1, \\
0, &  \mbox{ if } k>1,
\end{array} \right.
$$ with the usual convention of replacing $(\beta^{h})^k$ by $\beta_{k, q}^{h}$.

In this paper we construct $q$-Bernoulli numbers with weight
$\alpha$ which are different the extension of Carlitz's
$q$-Bernoulli numbers. By using $p$-adic $q$-integral equations on
$\mathbb{Z}_p$, we investigate some interesting identities and
relations on  the  $q$-Bernoulli numbers with weight $\alpha$.

\begin{center} {\bf 2. On the weighted  $q$-Bernoulli numbers and polynomials  } \end{center}

Let $f_n(x)=f(x+n)$. By (1), we see that

$$\aligned
 q I_q(f_1) & =   q \lim_{N \to \infty}\dfrac{1}{[p^N]_q}\sum_{x=0}^{ p^N-1}f(x+1)q^x,\\
 & = \lim_{N \to \infty}\dfrac{1}{[p^N]_q}\sum_{x=0}^{ p^N-1}f(x)q^x+ \lim_{N \to \infty}\dfrac{f(p^N)q^{p^N}-f(0)}{[p^N]_q}
 \\
 &= \int_{\mathbb{Z}_p} f(x) d\mu_{q}(x)  + (q-1)f(0)+
 \dfrac{q-1}{\log q}f^\prime (0),
\endaligned \eqno(3)
$$ and
$$\aligned
 q^2 I_q(f_2) & =   q^2  \int_{\mathbb{Z}_p} f(x+2) d\mu_{q}(x)  \\
 &=\lim_{N \to \infty}\dfrac{1}{[p^N]_q}\sum_{x=0}^{ p^N-1}f(x+1)q^{x+1}+ \lim_{N \to
 \infty}\dfrac{f(1+p^N)q^{p^N+1}-f(1)q }{[p^N]_q}\\
 & = q \int_{\mathbb{Z}_p} f(x+1) d\mu_{q}(x)  + q (q-1)f(1)+ q
 \dfrac{q-1}{\log q}f^\prime (1) \\
 &= I_q(f)+ (q-1)f(0) + q (q-1)f(1)+ \dfrac{q-1}{\log q}f^\prime (0)+
 \dfrac{q(q-1)}{\log q}f^\prime (1) . \endaligned
$$
 Thus,  we have
$$q^2 I_q(f_2)-I_q(f)= \sum_{l=0}^1 q^l (q-1) f(l)+ \dfrac{q-1}{\log q}\sum_{l=0}^1 q^l  f^\prime(l) .$$
Continuing this process, we obtain the following theorem.

\bigskip
{ \bf Theorem 1.} For $ n \in \Bbb N $, we have
$$q^n I_q(f_n)-I_q(f)= (q-1) \sum_{l=0}^{n-1} q^l f(l)+ \dfrac{q-1}{\log q}\sum_{l=0}^{n-1} q^l  f^\prime(l) .$$
In particular, $n=1$,
$$q I_q(f_1)-I_q(f)= (q-1) f(0) + \dfrac{q-1}{\log q}f^\prime (0) .$$
\bigskip

For $\alpha \in \Bbb N$, we evaluate  the following bosonic
$p$-adic $q$-integral on $\Bbb Z_p$:

$$  \widetilde{\beta}_{n, q}^{(\alpha)}=\int_{\mathbb{Z}_p} [x]_{q^\alpha}^n d\mu_{q}(x) =
\dfrac{1-q}{(1-q^\alpha)^n} \sum_{l=0}^n \binom nl (-1)^l \dfrac{
\alpha l+ 1}{1-q^{\alpha l+1}}. \eqno(4)$$

From (4), we note that

$$\aligned
 &\dfrac{1-q}{(1-q^\alpha)^n} \sum_{l=0}^n \binom nl (-1)^l \dfrac{
\alpha l+ 1}{1-q^{\alpha l+1}}  \\
 &=\dfrac{1-q}{(1-q^\alpha)^n} \sum_{l=0}^n \binom nl (-1)^l \dfrac{
\alpha l}{1-q^{\alpha l+1}} + \dfrac{1-q}{(1-q^\alpha)^n}
\sum_{l=0}^n \binom nl (-1)^l \dfrac{1}{1-q^{\alpha l+1}}  \\
 & =\dfrac{ n \alpha (1-q)}{(1-q^\alpha)^n} \sum_{l=1}^{n} \binom{ n-1}{l-1}  \dfrac{(-1)^l
}{1-q^{\alpha l+1}} + \dfrac{1-q}{(1-q^\alpha)^n}
\sum_{m=0}^\infty q^m (1-q^{m \alpha})^n \\
&=-n \alpha \dfrac{ (1-q)}{(1-q^\alpha)^n} \sum_{m=0}^\infty q^{m
\alpha+ m}(1- q^{m \alpha} )^{n-1} + (1-q) \sum_{m=0}^\infty q^m
[m]_{q^{\alpha}}^n\\
&= \dfrac{ - n \alpha }{[\alpha]_q } \sum_{m=0}^\infty q^{m
\alpha+ m}[m]_{q^{\alpha}}^{n-1}+ (1-q) \sum_{m=0}^\infty q^m
[m]_{q^{\alpha}}^n.
\endaligned
\eqno(5)
$$
Therefore, by (4) and (5), we obtain the following theorem.

\bigskip
{ \bf Theorem 2.} Let $ n, \alpha  \in \Bbb N $. Then we  have
$$-\dfrac{ \widetilde{\beta}_{n, q}^{(\alpha)}}{n}=\dfrac{  \alpha }{[\alpha]_q } \sum_{m=0}^\infty q^{m
\alpha+ m}[m]_{q^{\alpha}}^{n-1}- \dfrac{(1-q)}{n}
\sum_{m=0}^\infty q^m [m]_{q^{\alpha}}^n.$$
\bigskip

For $\alpha =1$, we note that $ \widetilde{\beta}_{n, q}^{(1)}$
are  same the Carlitz $q$-Bernoulli numbers.  In this paper, $
\widetilde{\beta}_{n, q}^{(\alpha)}$ are called the $n$-th
$q$-Bernoulli numbers with weight $\alpha$.

By (4) and (5), we easily get

$$ \int_{\mathbb{Z}_p} e^{[x]_{q^\alpha} t} d\mu_{q}(x)
= -t \dfrac{  \alpha }{[\alpha]_q } \sum_{m=0}^\infty q^{m \alpha+
m} e^{[m]_{q^{\alpha}} t}+(1-q) \sum_{m=0}^\infty q^m
e^{[m]_{q^{\alpha}}t}. \eqno(6)$$

From (6), we have

$$\sum_{n=0}^\infty \widetilde{\beta}_{n, q}^{(\alpha)} \dfrac{t^n}{n!} =-t \dfrac{  \alpha }{[\alpha]_q } \sum_{m=0}^\infty q^{m \alpha+
m} e^{[m]_{q^{\alpha}} t}+(1-q) \sum_{m=0}^\infty q^m
e^{[m]_{q^{\alpha}}t}.$$

Therefore, we obtain the following corollary.

\bigskip
{ \bf Corollary 3.} Let $ F_q^{(\alpha)}(t)= \sum_{n=0}^\infty
\widetilde{\beta}_{n, q}^{(\alpha)} \dfrac{t^n}{n!}.$ Then we have
$$F_q^{(\alpha)}(t)=-t \dfrac{  \alpha }{[\alpha]_q } \sum_{m=0}^\infty q^{m \alpha+
m} e^{[m]_{q^{\alpha}} t}+(1-q) \sum_{m=0}^\infty q^m
e^{[m]_{q^{\alpha}}t}.$$
\bigskip

 Now, we consider the $q$-Bernoulli
polynomials with weight $\alpha$ as follows:

$$ \widetilde{\beta}_{n, q}^{(\alpha)}(x)  = \int_{\mathbb{Z}_p} [x+y]_{q^\alpha}^n
d\mu_{q}(y),
 \text{ for } n \in \Bbb Z_+ \text{ and } \alpha \in \Bbb N.
 \eqno(7)
$$

By (7), we see that

$$\aligned
  \widetilde{\beta}_{n, q}^{(\alpha)}(x)  &= \dfrac{1-q}{(1-q^\alpha)^n} \sum_{l=0}^n \binom nl (-1)^l
q^{\alpha l x}  \dfrac{
\alpha l+ 1}{1-q^{\alpha l+1}} \\
&= - n \dfrac{  \alpha }{[\alpha]_q } \sum_{m=0}^\infty q^{m
\alpha+ m}[m+x]_{q^{\alpha}}^{n-1}+ (1-q) \sum_{m=0}^\infty q^m
[m+x]_{q^{\alpha}}^n.
\endaligned
\eqno(8)
$$
Let $ F_q^{(\alpha)}(t,x)= \sum_{n=0}^\infty \widetilde{\beta}_{n,
q}^{(\alpha)}(x) \dfrac{t^n}{n!}.$
Then we have
$$\aligned
F_q^{(\alpha)}(t, x)& = -t \dfrac{  \alpha }{[\alpha]_q }
\sum_{m=0}^\infty q^{m \alpha+ m} e^{[m+x]_{q^{\alpha}} t}+(1-q)
\sum_{m=0}^\infty q^m e^{[m+x]_{q^{\alpha}}t}\\
&= \sum_{n=0}^\infty \widetilde{\beta}_{n, q}^{(\alpha)}(x)
\dfrac{t^n}{n!}.\endaligned \eqno(9)$$ By simple calculation, we
easily get

$$ \aligned
 \widetilde{\beta}_{n, q}^{(\alpha)}(x) &= \sum_{l=0}^n \binom nl
[x]_{q^\alpha}^{n-l} q^{ \alpha l x} \int_{\mathbb{Z}_p}
[y]_{q^\alpha}^l d\mu_{q}(y)\\
&= \sum_{l=0}^n \binom nl [x]_{q^\alpha}^{n-l} q^{ \alpha l x}
 \widetilde{\beta}_{n, q}^{(\alpha)}.
\endaligned
\eqno(10)
$$
Therefore, by (8), (9) and (10), we obtain the following theorem.

\bigskip
{ \bf Theorem 4.} For $ n \in \Bbb Z_+ $ and  $ \alpha \in \Bbb
N$, we have
$$\aligned
\widetilde{\beta}_{n, q}^{(\alpha)}(x) &=
\dfrac{(1-q)}{(1-q^\alpha)^n} \sum_{l=0}^n \binom nl (-1)^l
q^{\alpha l x}  \dfrac{
\alpha l+ 1}{1-q^{\alpha l+1}} \\
&= - n \dfrac{  \alpha }{[\alpha]_q } \sum_{m=0}^\infty q^{m
\alpha+ m}[m+x]_{q^{\alpha}}^{n-1}+ (1-q) \sum_{m=0}^\infty q^m
[m+x]_{q^{\alpha}}^n.
\endaligned
$$
Moreover,
$$\widetilde{\beta}_{n, q}^{(\alpha)}(x) =\sum_{l=0}^n \binom nl [x]_{q^\alpha}^{n-l} q^{ \alpha l x}
\widetilde{\beta}_{l, q}^{(\alpha)}.$$

\bigskip

By Theorem 1, we see that
$$q^n \widetilde{\beta}_{m, q}^{(\alpha)}(n) -\widetilde{\beta}_{m, q}^{(\alpha)}
 = (q-1) \sum_{l=0}^{n-1} q^l [l]_{q^{\alpha}}^m +
  m \dfrac{ \alpha }{[\alpha]_q}\sum_{l=0}^{n-1} q^{\alpha l+l} [l]_{q^{\alpha}}^{m-1} .$$

Therefore, we obtain the following theorem.

\bigskip
{ \bf Theorem 5.} For $ m \in \Bbb Z_+ $ and  $ \alpha, n \in \Bbb
N $, we have $$ q^n \widetilde{\beta}_{m, q}^{(\alpha)}(n)
-\widetilde{\beta}_{m, q}^{(\alpha)}
 = (q-1) \sum_{l=0}^{n-1} q^l [l]_{q^{\alpha}}^m +
 \dfrac{ m \alpha }{[\alpha]_q}\sum_{l=0}^{n-1} q^{\alpha l+l} [l]_{q^{\alpha}}^{m-1} .$$

\bigskip

In (3), it is known that
$$ q I_q(f_1)- I_q(f)= (q-1)f(0)+ \dfrac{q-1}{\log q}f^\prime
(0).$$

If we take $f(x)=e^{[x]_{q^{\alpha}} t}$, then we have
$$\aligned
  (q-1)+ \dfrac{\alpha}{[\alpha]_q}  t &
  = q \int_{\mathbb{Z}_p} e^{ [x+1]_{q^\alpha} t} d\mu_{q}(x) -  \int_{\mathbb{Z}_p} e^{ [x]_{q^\alpha} t}
 d\mu_{q}(x)\\
 &= \sum_{n=0}^\infty \left(
  q \widetilde{\beta}_{n, q}^{(\alpha)}(1) - \widetilde{\beta}_{n, q}^{(\alpha)} \right)
  \dfrac{t^n}{n!}.
\endaligned
\eqno(11)
$$
Therefore, by (11), we obtain the  following theorem.

\bigskip
{ \bf Theorem 6.} For $ \alpha \in \Bbb N$ and $ n \in \Bbb Z_+ $,
we have
$$\widetilde{\beta}_{0, q}^{(\alpha)}=1, \text{  and  } \quad  q \widetilde{\beta}_{n, q}^{(\alpha)}(1) - \widetilde{\beta}_{n, q}^{(\alpha)}  =
\left \{\begin{array}{ll}
\dfrac{\alpha}{[\alpha]_q}, & \mbox{ if } n=1, \\
0, &  \mbox{ if } n>1.
\end{array} \right.
$$
\bigskip

By (10) and Theorem 6,  we obtain the following corollary.

\bigskip
{ \bf Corollary 7.}  For $ \alpha \in \Bbb N$ and $ n \in \Bbb Z_+
$, we have
$$\widetilde{\beta}_{0, q}^{(\alpha)}=1, \text{  and  } \quad
 q ( q^{\alpha} \widetilde{\beta}^{(\alpha)}+1)^n - \widetilde{\beta}_{n, q}^{(\alpha)}  =
\left \{\begin{array}{ll} \dfrac{\alpha}{[\alpha]_q}, & \mbox{ if
} n=1, \\ \\
0, &  \mbox{ if } n>1,
\end{array} \right.
$$
 with the usual convention of replacing $ (\widetilde{\beta}_q^{(\alpha)})^n $
by $  \widetilde{\beta}_{n, q}^{(\alpha)}$.

\bigskip

From (7), we can easily derive the following equation (12).

$$\aligned
  \int_{\mathbb{Z}_p}  [x+y]_{q^\alpha}^n d\mu_{q}(y) & = \dfrac{[d]_{q^\alpha}^n}{[d]_q}
  \sum_{a=0}^{d-1} q^a  \int_{\mathbb{Z}_p}  \left[ \dfrac{x+a}{d}+ y \right]_{q^{\alpha
  d}}^n
 d\mu_{q^{d}}(y)\\
 &= \dfrac{[d]_{q^\alpha}^n}{[d]_q}
  \sum_{a=0}^{d-1} q^a
   \widetilde{\beta}_{n, q^d}^{(\alpha)} \left( \dfrac{x+a}{d} \right). \endaligned
\eqno(12)
$$
Therefore, by (12), we obtain the following theorem.

\bigskip
{ \bf Theorem 8.} For $ n \in \Bbb Z_+ $ and  $ \alpha, d \in \Bbb
N$, we have
$$\widetilde{\beta}_{n, q}^{(\alpha)}(x)= \dfrac{[d]_{q^\alpha}^n}{[d]_q}
  \sum_{a=0}^{d-1} q^a
   \widetilde{\beta}_{n, q^d}^{(\alpha)} \left( \dfrac{x+a}{d} \right).
$$
\bigskip

From (7), we note that

$$\aligned
  \widetilde{\beta}_{n, q^{-1}}^{(\alpha)}(1-x) &= \int_{\mathbb{Z}_p}  \left[1-x+x_1\right]_{q^{-\alpha}}^n d\mu_{q^{-1}}(x_1)\\
  & = (-1)^n  q^{ \alpha n}\int_{\mathbb{Z}_p}  [x+x_1]_{q^{\alpha}}^n d\mu_{q}(x_1) \\
 &= (-1)^n  q^{ \alpha n} \widetilde{\beta}_{n, q}^{(\alpha)}(x). \endaligned
\eqno(13)
$$

Therefore, by (13), we obtain the following theorem.

\bigskip
{ \bf Theorem 9.} For $ n \in \Bbb Z_+$ and $ \alpha \in  \Bbb N$,
we have
$$
  \widetilde{\beta}_{n, q^{-1}}^{(\alpha)}(1-x) = (-1)^n  q^{ \alpha n} \widetilde{\beta}_{n, q}^{(\alpha)}(x).
$$
\bigskip

It is easy to show that

$$\aligned
  \int_{\mathbb{Z}_p}  [1-x ]_{q^{-\alpha}}^n d\mu_{q}(x) &=
   \int_{\mathbb{Z}_p}  (1-[x]_{q^{\alpha}})^n d\mu_{q}(x) \\
 &= (-1)^n  q^{ \alpha n} \int_{\mathbb{Z}_p}  [x-1 ]_{q^{\alpha}}^n d\mu_{q}(x). \endaligned
\eqno(14)
$$

By (13) and (14), we obtain the following corollary.

\bigskip
{ \bf Corollary 10.} For $ \alpha \in \Bbb N $ and $  n \in \Bbb
Z_+$, we have
$$\aligned \int_{\mathbb{Z}_p}  [1-x ]_{q^{-\alpha}}^n d\mu_{q}(x) & = \sum_{l=0}^n \binom nl (-1)^l
 \widetilde{\beta}_{l, q}^{(\alpha)} \\
 &= (-1)^n  q^{ \alpha n} \widetilde{\beta}_{n, q}^{(\alpha)}(-1) \\
 &=   \widetilde{\beta}_{n, q^{-1}}^{(\alpha)}(2).\endaligned
$$
\bigskip

From  Theorem 4, we have

$$\aligned
  & q^2 \widetilde{\beta}_{n, q}^{(\alpha)}(2)- n q^{1+\alpha} \dfrac{ \alpha}{[\alpha]_q}-q^2+q \\
 &= q^2  \sum_{l=0}^n \binom nl q^{ \alpha l}
  \widetilde{\beta}_{l, q}^{(\alpha)}(1)- n q^{1+\alpha} \dfrac{ \alpha}{[\alpha]_q}-q^2+q  \\
 &=q^2 + q^2 n q^{\alpha}\widetilde{\beta}_{1, q}^{(\alpha)}(1)+q^2 \sum_{l=2}^n \binom nl q^{ \alpha l}
  \widetilde{\beta}_{l, q}^{(\alpha)}(1) - n q^{1+\alpha} \dfrac{ \alpha}{[\alpha]_q}-q^2+q  \\
  &=q^2 + q^{1+\alpha} n  \left( \dfrac{ \alpha}{[\alpha]_q}+ \widetilde{\beta}_{1, q}^{(\alpha)}\right)
  +q  \sum_{l=2}^n \binom nl q^{ \alpha l}
  \widetilde{\beta}_{l, q}^{(\alpha)} - n q^{1+\alpha} \dfrac{ \alpha}{[\alpha]_q}-q^2+q  \\
 &= q  \sum_{l=0}^n \binom nl q^{ \alpha l}  \widetilde{\beta}_{l, q}^{(\alpha)}
 = q ( q  \widetilde{\beta}^{(\alpha)} +1)^n = \widetilde{\beta}_{n, q}^{(\alpha)} \text{ if } n>1.\endaligned
$$

Therefore, we obtain the following theorem.

\bigskip
{ \bf Theorem 11.} For $ n \in \Bbb Z_+$  with  $ n >1$, we have
$$
  q^2 \widetilde{\beta}_{n, q}^{(\alpha)}(2)
  =  n q^{1+\alpha} \dfrac{ \alpha}{[\alpha]_q}+q^2-q + \widetilde{\beta}_{n, q}^{(\alpha)}.
$$
\bigskip

\begin{center}{\bf REFERENCES}\end{center}

\begin{enumerate}

\item
{ L.  Carlitz}, { Expansion of $q$-Bernoulli numbers and
polynomials},
 { Duke Math. J.} { 25}(1958), 355-364.

\bibitem{}
{ L.  Carlitz}, { $q$-Bernoulli numbers and polynomials},
 { Duke Math. J.} {15}(1948), 987-1000.

\item
{ T. Kim,}  {   $q$-Bernoulli numbers and polynomials associated
with Gaussian binomial coefficients},{ Russ. J. Math. Phys.}
{15}(2007), 51-57.

\item
{T. Kim,}  {  Some identities on the $q$-Euler polynomials of
higher order and $q$-Stirling numbers by the fermionic $p$-adic
integral on $\mathbb{Z}_p$,} { Russ. J.  Math. phys.}{ 16}(2009),
484-491.

\item
{T. Kim,} { Barnes type multiple $q$-zeta function and $q$-Euler
polynomials,} { J. phys. A : Math.  Theor.} {43}(2010) 255201,  11
pages.

\item
{T. Kim,}  {  Note on the Euler $q$-zeta functions,} {J.  Number
Theory} { 129}(2009), 1798-1804.

\item
{T. Kim,}  {  On a $q$-analogue of the $p$-adic log gamma
functions and related integrals,} {J. Number Theory} {76}(1999),
320-329.

\item
{  L. C. Jang, W.-J. Kim,  Y. Simsek}, { A study on the $p$-adic
integral representation on $\Bbb Z_p$ associated with Bernstein
and Bernoulli  polynomials,}  { Advances in Difference Equations},
2010(2010), Article ID 163217, 6 pages.

\item
{  Y. Simsek, M. Acikgoz}, { A new generating function of
$q$-Bernstein-type polynomials and their interpolation function,}
{ Abstract and Applied Analysis}, 2010(2010), Article ID 769095,
12 pages.

\item
{ H. Ozden, I. N. Cangul,  Y. Simsek}, { Remarks on $q$-Bernoulli
numbers associated with Daehee numbers,} { Adv. Stud. Contemp.
Math.} {18}(2009), 41-48.

\item
{ A. S. Hegazi, M. Mansour}, { A note on $q$-Bernoulli numbers and
polynomials,} { J. Nonlinear  Math. Phys.} 13(2006), 9-18.

\item
{   M. Acikgoz, Y. Simsek}, { On multiple interpolation function
of the N\"{o}rlund-type $q$-Euler polynomials,} { Abst. Appl.
Anal.} Article ID 382574(2009), 14 pages.

\item
{H. Ozden, Y. Simsek,}  { A new extension of $q$-Euler numbers and
polynomials related to their interpolation functions,} { Appl.
Math. Letters} {21}(2008), 934-938.

\item
{L. C. Jang,} {A  note on N\"{o}rlund-type  twisted $q$-Euler
polynomials and numbers of higher order associated with fermionic
invariant $q$-integrals,} { J. Inequal. Appl.} 2010(2010), Article
ID 417452, 12 pages.

\item
{S.-H. Rim, J.-H. Jin, E.-J. Moon,  S.-J. Lee}, { On multiple
interpolation function of the $q$-Genocchi  polynomials}, {J.
Inequal. Appl.} 2010(2010), Article ID 351419, 13 pages.

\item
{  B. A. Kupershmidt}, { Reflection symmetries  of $q$-Bernoulli
polynomials,} { J. Nonlinear  Math. Phys.} 12(2005), 412-422.

\end{enumerate}

\end{document}